\documentstyle[12pt]{article}

\def\ni{\noindent}
\def\qed{\ \vrule width2mm height3mm \vspace{3mm}}

\def\ss #1 {\raisebox{-1.0ex}{$#1$}}
\def\std {Gr\"{o}bner }

\def\Proj{\mbox{\bf P}}
\def\arrow{\rightarrow}

\newenvironment{define}{\vspace{3mm}\ni{\bf Definition:}}{\vspace{3mm}}

\newenvironment{rem}{\vspace{3mm}\ni{\bf Remark:}}{\vspace{3mm}}

\newenvironment{note}{\ni {\em Note:}}{}

\newtheorem{thm}{Theorem}
\newtheorem{lem}[thm]{Lemma}
\newtheorem{prop}[thm]{Proposition}
\newtheorem{cor}[thm]{Corollary}

\newenvironment{ex}{\vspace{3mm}\ni{\bf Example}}{\vspace{3mm}}
\def\P #1 { {\bf P}^#1 }
\def\prs{\P {r_1} \times \cdots \times \P {r_s}}
\def\to {\rightarrow}

\begin{document}

\title{ Initial ideals, Veronese subrings, and rates of algebras}
\author{ (To appear in Advances in Mathematics)\\
David Eisenbud \thanks{The authors are grateful to the
NSF for support during the preparation of this paper.  The second and third
authors were
supported by NSF Postdoctoral Fellowships.}
\\ Dept. of Math., Brandeis Univ., Waltham MA 02254 \\
eisenbud@math.brandeis.edu
   \and Alyson Reeves
\\ Dept. of Math., Brandeis Univ., Waltham MA 02254 \\
reeves@math.brandeis.edu
	\and Burt Totaro \\ Dept. of Math., Univ. of Chicago, Chicago IL
  60637 \\ totaro@math.uchicago.edu}

\maketitle

\begin{abstract}
 	Let $S$ be a polynomial ring over an infinite field and let $I$ be a
 homogeneous ideal of $S$.  Let $T_d$ be a polynomial ring
 whose variables correspond to the monomials of degree
 $d$ in $S$.  We study the initial ideals of the ideals $V_d(I) \subset
 T_d$ that define the Veronese subrings of $S/I$.
In suitable orders, they
 are easily deduced from the initial ideal
  of $I$.
 We show that
$in(V_d(I))$ is generated in degree $\leq max(\lceil reg(I)/d \rceil,2)$, where
$reg(I)$ is the regularity of the ideal $I$. (In other words, the $d^{th}$
Veronese subring of any commutative graded ring $S/I$ has a \std basis
of degree $\leq max(\lceil reg(I)/d \rceil, 2)$.)
We also give bounds on the regularity of $I$ in terms of the
degrees of the generators of $in(I)$ and some combinatorial
data.
This implies a
version of Backelin's Theorem that
 high Veronese subrings of any ring are homogeneous Koszul algebras
 in the sense of Priddy [Pr70].

We also give a general obstruction for a homogeneous ideal $I\subset S$ to have
an initial ideal $in(I)$ that is generated by quadrics, beyond the obvious
requirement that $I$ itself should be generated by quadrics, and the
stronger statement that $S/I$ is Koszul.
We use the obstruction to show that in certain dimensions,
a generic complete intersection of quadrics cannot have an initial ideal
that is generated by quadrics.

For the application to Backelin's Theorem, we require a result of
Backelin whose proof has never appeared.  We give a simple proof of
a sharpened version, bounding the rate of growth of the degrees of
generators for syzygies of any multihomogeneous module over a
polynomial ring modulo an ideal generated by monomials, following
a method of Bruns and Herzog.

\end{abstract}

  	Notation:  Throughout this paper we write $S = k[x_1,\ldots, x_r]$ for
 the graded polynomial ring in $r$ variables over an infinite field $k$.
 We will generally deal with a monomial order $>$ on $S$.  We always suppose
$x_1>\cdots >x_r$; by the initial term $in_>(p)$ of a polynomial $p$ we mean
the term with
 the largest monomial.  Similarly, by the {\bf initial ideal} $in_>(I)$ of $I$,
 we mean the ideal generated by the initial terms of all polynomials in $I$:
	$$in_>(I) = \langle in_>(p) \,| \, p \in I \rangle.$$
We write $T_d$ for the polynomial ring over $k$
 whose variables correspond to the monomials of degree $d$ in $S$; there
 is a natural map $\phi_d: T_d \rightarrow S$ sending each variable of $T_d$
 to the
 corresponding monomial in $S$.  If $I \subset S$ is an ideal, we write
$V_d(I)$ for
 the preimage of $I$ in $T_d$, and $A_{(d)}$ for $T_d/V_d(I)$.

Furthermore, we make the following definitions:

\begin{define}
For a homogeneous ideal $I \subset S$, the {\bf minimal generators} of $I$ are
the homogeneous elements of $I$ not in $(x_1,\ldots,x_r)I$.
Let {\bf $\delta(I)$}  be the maximum of the degrees of minimal generators
of $I$,
and let
{\bf $\Delta(I)$} be the minimum, over all choices of variables and of
monomial orderings of $S$ of the maximum of the degrees of minimal generators
 of the
initial ideal of $I$.
\end{define}

All ideals and rings that appear will be graded.

\section{ Introduction and motivating results}
\label{intro}

Given a homogeneous ideal $I \subset S$ it is a matter of both
computational and theoretical interest to know how low
the degree of $in_>(I)$ can be made by choosing variables
and monomial order on $S$ in an appropriate way.  In particular,
one may ask which ideals $I$ admit {\it quadratic}\ initial ideals;
that is, for which $I$ are there choices of variables and
order such that $in_>(I)$ can be generated by monomials of
degree 2?

One of the theoretical reasons for  interest in this question is
that if $I$ admits a quadratic initial ideal then, by a result
of Fr\"oberg and a deformation argument noticed by
Kempf and others, $A:=S/I$ is
a (homogeneous) Koszul algebra in the sense of Priddy [Pr70];
that is, the residue field $k$ of $A$ admits an $A$-free resolution
whose maps are given by matrices of linear forms.  Using complex
arguments about a lattice of ideals derived from a presentation
of $A$ as a quotient of a free noncommutative algebra, Backelin
[Ba86] proved that for any graded ring $A$ as above the
{\bf $d^{th}$ Veronese subring }
  $$	A_{(d)} := \bigoplus_{i=1}^{\infty} A_{di},$$
is Koszul for all sufficiently large $d$.  Our work started
from a request by George Kempf for a simpler proof.
In this paper we shall prove the
stronger result that $	A_{(d)}$ admits a quadratic initial
ideal for all sufficiently large $d$.

To make these results more quantitative, we define
a measure of the rate of growth of the degrees of the
syzygies in a minimal free resolution:

\begin{define}
Let
$A = k \oplus A_1 \oplus A_2 \oplus \cdots$ be a graded ring.
For any finitely generated graded $A$-module $M$,
set $t_i^A(M) = max\{j\ | \ Tor^A_i(k,M)_j \neq 0\}$, where
$Tor_i^A(k,M)_j$ denotes the $j^{th}$ graded piece of $Tor_i^A(k,M)$.
The {\bf rate} of $A$ is defined by Backelin [Ba86] to be

 $$rate(A) = sup \{ (t^A_i (k)-1)/(i-1)  |   i \geq 2  \}.$$
\end{define}

For example, $A$ is Koszul iff
$rate(A) = 1$.  It turns out that the rate of any graded algebra
is finite (see for example [An86]) and Backelin actually
proves

\begin{thm}
\label{thm0}
 ([Ba86]):  $rate(A_{(d)}) \leq max(1, rate(A)/d)$.
\end{thm}

One should compare this with the rather trivial result
(Proposition \ref{mum-prop} below)
that if a homogeneous ideal $I$ can be generated by
forms of degree $m$, then the ideal $V_d(I)$ defining
$A_{(d)}$ can be generated by forms of degree
$\leq max (2, \lceil m/d \rceil)$.
In the notation introduced above,
$\delta(V_d(I)) \leq max (2,\lceil \delta(I)/d \rceil)$.
A similar result with $\Delta$ (the minimum, over all choices
of variables and of monomial orderings $>$, of the maximum degree
of a minimal generator of $in_>(I)$) would lead to a bound on the
rate by virtue of Proposition \ref{prop4}.
Unfortunately,
as we show in Example \ref{counterexample} below, it is
{\it not} true that if some initial ideal of $I$ can be
generated by monomials of degree $m$ then $V_d(I)$ admits
an initial ideal generated by forms of degree
$\leq max (2, \lceil m/d \rceil)$. But there {\it is}
a replacement
for $m$ that makes such a formula true, and this is the main
result of this paper.

Recall that the {\bf Castelnuovo-Mumford regularity} of $I$
is defined as follows:
\begin{define}
For $I \subset S$,
 the {\bf regularity} of $I$ is defined as

 $$reg(I) = max \{t^S_i (I) - i  |  i \geq 0 \}.$$
\end{define}

Since
$t^S_0(I) = \delta(I)\leq reg(I)$, the regularity
is $\geq$ the maximal degree of the generators
of $I$.  One may think of the regularity as
a more stable measure of the size of the generators of $I$.  Our
main result is that we may replace the degree of the generators
of $I$ by the regularity and get a bound on the degrees of the
initial ideals of Veronese powers:

\begin{thm}
\label{reg-thm}
        $$\Delta(V_d(I)) \leq max (2, \lceil reg (I)/d \rceil). $$
 In particular, if $d \geq  reg( I)/2 $ then $\Delta(V_d(I)) = 2$.
\end{thm}
In section \ref{multi} we
explain how to generalize this result
to Segre products of Veronese embeddings.

To deduce a version of Backelin's Theorem, one needs a result
strengthening the theorem of Fr\"oberg mentioned above.  Such
a result was stated without proof by Backelin ([Ba86, pp98ff]):

\begin{prop}
\label{prop4}
 If $A = S/I$ with $I$ a homogeneous ideal, then
 $rate(A) \leq \Delta(I) - 1$.  In particular, if $\Delta(I) = 2$ then $A$
 is Koszul.
\end{prop}

We will give a simple proof of
this proposition (and something more general)
in section \ref{resolution}
following ideas of Bruns, Herzog, and Vetter [BrHeVe].

Unfortunately the converse of this result is not true:  in particular,
the algebra $A$ may be Koszul without $I$ admitting a
quadratic initial ideal. In section \ref{obstsec},
we formulate another obstruction for an ideal $I\subset S$ to have
an quadratic initial ideal. An easily stated part of Theorem \ref{obstruction}
is that if $I$ admits a quadratic initial ideal then $I$ contains
far more quadrics of low rank than would a generic subspace of quadrics.
We may make this quantitative as follows:
\begin{cor}
\label{cor to Lemma 18}
 If $I\subset S$ admits a quadratic initial ideal,
and $dim (S/I) = n$, then
$I$ contains an $m$-dimensional space of quadrics of rank
$\leq 2(n+m)-1$ for every $m \leq codim(I)$.
\end{cor}
The obstruction
shows that in certain dimensions, a generic complete intersection of quadrics
has no initial ideal $in(I)$ which is generated by quadrics, even though
every complete intersection of quadrics is a Koszul algebra.
There seems no reason to believe that this obstruction, even with
the Koszul condition, is enough to guarantee that an ideal admits a
quadratic initial ideal, so we pose as a problem the question
raised at the
beginning of the introduction :

\smallskip
{\sl Find necessary and sufficient conditions for an ideal to
admit a quadratic initial ideal.}
\smallskip

A noncommutative analogue of some of the results on initial
ideals, in which $T_d$
is replaced by a free noncommutative algebra mapping onto $S$,
is given in [Ei].

\section{Initial ideals for Veronese subrings}
\label{vero}

        As above, let
$$T_d = k[\{z_{m}\}] \mbox{ where $m$ is a monomial of $S$ of
degree $d$},$$
and let $\phi_d: T_d \rightarrow S$ be the map sending $z_{m}$ to
$m$.
  If $J \subset S$ is a
 homogeneous ideal, let $V_d(J)$ denote the preimage of $J$ in $T_d$.  It is
 easy to see that $V_d(J)$ is generated by the kernel of $\phi_d$ and, for each
 generator $g$ of $J$ in degree $e$, the preimages of the elements of degree
 $nd$ in $(x_1,\ldots,x_r)^{nd-e}g$, where $nd$ is the smallest multiple of $n$
that is
 $\geq e$.  These elements have degree $n$ in $T_d$.  Since ${\rm ker}(\phi_d)$
 is
generated
 by forms of degree 2 it follows that $V_d(J)$ is generated by forms of
 degree $\leq max(\lceil \delta(J)/d \rceil,2)$.

This gives a proof of the following well-known Proposition.

\begin{prop}
\label{mum-prop}
 $$\delta(V_d(I)) \leq max(2, \lceil \delta(I)/d \rceil).$$ In particular,
if $d \geq  \delta(I)/2 $
 then $V_d(I)$ is generated by quadrics.
\end{prop}

This Proposition is mentioned by Mumford in [Mu70] (in a slightly
different form), though it is surely much older.

        We extend the given monomial order on $S$ to a monomial
 order on $T_d$ as follows:   If $a$, $b$ are monomials in $T_d$, then $a > b$
if
 $\phi(a) > \phi(b)$ or $\phi(a) = \phi(b)$ but $a$ is bigger than $b$ in the
reverse
 lexicographic order: that is, given two monomials in $T_d$ of the same degree
 having the same image in $S$ we order the factors of each in
 decreasing order, and take as larger the monomial with the smaller
 factor in the last place where the two differ. Here the order of the
 variables  $z_{m}$ is defined to be the same as the order in $S$ on the
monomials $m$.
We first compute the initial ideal
of ${\rm ker}(\phi_d)$.

\begin{prop}
\label{prop8}
   With notation as above, $in({\rm ker}(\phi_d)) \subset T_d$ is
 generated by quadratic forms for every $d$.
\end{prop}

\begin{rem}
 Barcanescu and Manolache [BaMa82] proved that the
 Veronese rings are Koszul -- which is a corollary of Proposition \ref{prop8}.
\end{rem}

 {\bf Proof}:  The ideal ${\rm ker}(\phi_d)$ is generated by quadratic forms,
 each a
 difference of two monomials that go to the same monomial under
 $\phi_d$.  Let $J$ be the monomial ideal generated by the initial terms of
 these quadratic elements of ${\rm ker}(\phi_d)$.  We have $J \subset
 in({\rm ker}(\phi_d))$,
and we
 wish to prove equality.  We will show that distinct monomials of $T_d$
 not in $J$ map by $\phi_d$ to distinct monomials of $S$.  It will follow that,
 for each $e$,
\begin{eqnarray*}
        dim S_{de} &\geq &dim (T_d/J)_e \cr
                   &\geq & dim (T_d/(in ({\rm ker}( \phi_d))))_e \cr
                  &=&dim (T_d/({\rm ker}( \phi_d)))_e \cr
                   &=& dim( S_{de}).
\end{eqnarray*}
  Thus $dim (T_d/J)_e = dim (T_d/(in( {\rm ker}( \phi_d))))_e$, and
$J = in({\rm ker}( \phi_d))$ as desired.

        Call the monomials not in $J$ ``standard", and say that a product
 of monomials of degree $d$ in $S$ is standard if its factors correspond to
 the factors of a standard monomial in $T_d$.
Since $J \subset in({\rm ker}(\phi_d))$, any monomial of $S_{de}$ may be
 written as
a standard product of $e$ monomials of degree $d$.
We must show that if $m \in  S$
 is a monomial of degree $de$, then there is a unique way of
 writing $m$ as a standard product $m_1\cdots  m_e$ of monomials of degree
$d$.
 We claim that this unique product is obtained by writing out the $de$
 factors of $m$ in decreasing order
        $$m = x_1\cdots x_1 x_2\cdots x_2 \cdots
x_r\cdots  x_r,$$
 and taking $m_1$ to be the product of the first $d$ factors,
$m_2$ to be the product of the next $d$ factors, and so on.

        First we prove the claim for a standard product $m_1m_2$ with
 just two factors. Suppose the sequences of indices of the two factors
 are
 $$i_1= (i_{11}\leq \cdots \leq i_{1d}), \;
 i_2= ( i_{21}\leq \cdots \leq i_{2d})$$
 and $m_1 > m_2$.  We must show that $ i_{1d} \leq  i_{21}$.
The product of the monomials $m_1'$ and $m_2
'$
obtained from $m_1$ and $m_2$ by interchanging the factors $x_{ i_{1d}}$
and $x_{ i_{21}}$
represents the same element of $S$.  The difference of these
 products represents a quadratic element of ${\rm ker}(\phi_d)$ . If
 $ i_{1d} >  i_{21}$ then $x_{ i_{1d}} < x_{ i_{21}}$
and thus $m_2' < m_2$.  Consequently
 the leading term of this quadratic element of ${\rm ker}(\phi_d)$ would
 correspond to $m_1m_2$, and the product $m_1m_2$ would not be standard.

        Now suppose that $m_1\cdot \cdots \cdot m_e$  is a standard product of
degree $de$,
 with $e$ arbitrary, and $m_1 \geq \cdots \geq m_e$.
Let $i_{j1}\leq \cdots \leq \alpha_{jd}$ be the indices of
 the variables in $m_j$.  If $i_{jd}>i_{j+1,1}$ for some $j$
then the product
 $m_jm_{j+1}$ would not be standard, contradicting the standardness of
 the entire product. \qed

\begin{note}
Even when $r = 3$ and $d = 2$, the initial terms of the
minimal system of generators for the ideal
${\rm ker}(\phi_2)$ do not generate $in_>({\rm ker}(\phi_2))$ under all
possible orders $>$.
In this case, the ideal is
 minimally generated by the $2 \times 2$ minors of the symmetric $3 \times 3$
matrix with entries $z_{ij}$, $i \leq j$.  There are 29
different initial ideals (depending upon the order chosen),
23 of which are generated entirely in degree 2. The other 6 each require an
additional generator of degree three.

It would be interesting to characterize the orders $>$ for which the ideal
$ in_>({\rm ker} (\phi_d))$ is generated by elements of degree 2.
\end{note}

 Following the proof of Proposition \ref{prop8} we say a monomial
 of $T_d$ is standard if it is not in $in({\rm ker} (\phi_d))$.
Let $\sigma: S \rightarrow T_d$ be
the $k$-linear map
 that takes each monomial to its unique standard representative.
  Because of the way we have defined the order on $T_d$ we have
 $in(\sigma(p)) = \sigma(in(p))$.  Also, $\sigma$ takes $J$ into $V_d(J)$.
As a consequence we have:
\begin{lem}
\label{lem9}
If $J$ is a homogeneous ideal of $S$,  then $in(V_d(J))$ is the
 ideal $K$ generated by $in({\rm ker} (\phi_d))$ and the monomials $\sigma(m)$
 for
$m \in in(J) \cap (im( \phi_d))$.
\end{lem}

{\bf Proof}:  For each degree $e$ it is clear that $dim_k(T_d/K)_e  \leq
 dim_k(S/J)_{de}$, so it is enough to show that $K \subset in(V_d(J))$.  Let
 $p \in J$
 be
 a form with initial term $m$.  Clearly $\sigma(p) \in V_d(J)$, and
$in(\sigma(p)) = \sigma(m)$.
\qed

\begin{define}
If $m \in S$ is a monomial, then
         following Eliahou and Kervaire [ElKe90] we
write $max(m)$ for the largest index $i$ such that $x_i$ divides $m$. We
call a monomial ideal $I$
{\bf (combinatorially) stable} if for every monomial $m \in I$ and
$j < max(m)$, the monomial
 $(x_j/x_{max(m)})m \in I$.
\end{define}

\begin{thm}
\label{char0-thm}
 If $J$ is a homogeneous ideal of $S$ such that $in(J)$ is
 combinatorially stable,  then $in(V_d(J))$ is
generated by $in({\rm ker}(\phi_d))$ and the monomials
 $\sigma(m)$ where $m$ runs over the minimal generators of $in(J) \cap
im (\phi_d)$.
  Thus if $\delta(in(J)) \leq u$, then
$\delta(in(V_d(J))) \leq max(2, \lceil u/d \rceil)$.
\end{thm}

{\bf Proof}:  Let $n \in in(J) \cap im (\phi_d)$.  By Lemma \ref{lem9}
it suffices to show that
 $\sigma(n)$ is divisible by some $\sigma(m)$, where $m$ is a minimal generator
of
 $in(J) \cap im (\phi_d)$.  Let $m'$ be an element of $in(J) \cap im (\phi_d)$
of minimal
 degree among those dividing $n$, and write $n = x_{i_1}\cdot \cdots \cdot
x_{i_{sd}}$
 with $i_1 \leq \cdots
 \leq i_{sd}$.  Say $deg \, m' = de$.
Since $in(J)$ is combinatorially stable, it follows that $m
 := x_{i_1}\cdot \cdots \cdot x_{i_{de}} \in in(J)$, and
since $m$ has degree $de$, we have $m \in
 im(\phi_d)$ as well. If $m$ were not a minimal generator of $in(J) \cap
im (\phi_d)$,
 then some proper divisor of it would be in $in(J) \cap (im (\phi_d))$ and
 would
divide $n$, contradicting our choice of $m'$.  As $\sigma(m)$ divides
$\sigma(n)$, we
 are done.
\qed

\begin{ex}
1: The hypothesis of stability cannot be dropped.  For example, if $r=3$ and
$J = (x_1x_2x_3)$, and we take the lexicographic or reverse lexicographic order
on $S$, then the initial ideal $in(V_d(J))$ (defined using the order on $T_d$
we associate to the given order on $S$)
 requires a cubic generator for all $d$.
\end{ex}

To obtain Theorem \ref{reg-thm} as given in the introduction, we need
to recall the notion of {\bf Castelnuovo-Mumford regularity}.

\begin{define}
For $I \subset S$,
 the {\bf regularity} of $I$ is defined as

 $$reg(I) = max \{t^S_i (I) - i  |  i \geq 0 \}.$$
\end{define}

\begin{note}
$t^S_0(I) = \delta(I)\leq reg(I)$.
\end{note}

Bayer and Stillman [BaSt87] give the following criterion for an ideal to
be $m$-regular, assuming (as we do throughout this paper) that the field
$k$ is infinite.

\begin{thm} [BaSt87]
\label{m-reg}
Let $I \subset S$ be an ideal generated in degrees $\leq e$.
The following conditions are equivalent:

\begin{enumerate}
\item $I$ is $e$-regular,

\item \begin{enumerate}
    \item For some $j\geq 0$ and for some linear forms
 $h_1,\ldots,h_j \in S_1$ we have
 $$((I,h_1,\ldots,h_{i-1}):h_i)_e = (I,h_1,\ldots,h_{i-1})_e $$
for $i=1,\ldots,j$, and
  \item $$(I,h_1,\ldots,h_j)_e = S_e.$$
  \end{enumerate}

\item Conditions 2a) and 2b) hold for some $j\geq 0$ and for {\em generic}
linear forms $h_1,\ldots ,h_j\in S_1$.

\end{enumerate}
\end{thm}

\begin{note}
Let $g$ be generically chosen in the Borel group $B$, the
subgroup of $Gl(r)$ consisting
of the upper triangular matrices.  Then $\langle gx_r,\ldots,gx_{r-j+1}
 \rangle$
is a generic linear subspace for $I$.  Since $gx_i$ is a generic linear
form with respect to $(I,gx_r,\ldots,gx_{i+1})$,
$x_i$ is a generic linear form with respect to $(g^{-1}I,x_r,\ldots,x_{i+1})$.
If $I$ is Borel-fixed, then $g^{-1}I = I$, and hence
we can replace $h_1,\ldots,h_j$ by
$x_r,\ldots,x_{r-j+1}$ in the statement of Theorem $\ref{m-reg}$ in this
case.
\end{note}

\begin{prop}
\label{reg-stab}
Let $I \subset S$ be a Borel-fixed monomial ideal generated in degrees
$\leq e$. Then
$I$ is $e$-regular if and only if $I_e$ is combinatorially stable.
\end{prop}

{\bf Proof:}
By the note following the statement of Theorem \ref{m-reg}, we may
replace $h_1,\ldots,h_j$ by $x_r,\ldots,x_{r-j+1}$ in the statement of the
theorem.

Suppose $I$ is $e$-regular.
Then 2b) of Theorem \ref{m-reg} implies
that $I$ includes all monomials in $x_1,\ldots,x_{r-j}$ of degree $e$.
And 2a) of the same theorem implies that for every monomial $m \in
I$ of degree $e$ with $max(m)>r-j$, $I$ also contains
$x_k/x_{max(m)} \cdot m$ for every $k$ with $1 \leq k \leq max(m)$.
Taken together,
these two statements imply $I_e$ is combinatorially stable.

Conversely, suppose $I_e$ is combinatorially stable,
and let $j$ be the smallest integer such that $I$ contains a power of
$x_{r-j}$. It follows that $x_{r-j}^e \in I_e$, and by stability,
$(x_1,\ldots,x_{r-j})^e \subset I$, and hence 2b)
holds. Let $m \in ((I,x_r,\ldots,x_{i+1}):x_{i})_{e}$ for some
$r-j+1 \leq i \leq r$; since $I$ is a monomial ideal, we can assume that
$m$ is a monomial in proving 2a).
If $m$ is divisible by $x_k$ for some
$i+1\leq k \leq r$, then it is clear that $m \in (I,x_r,\ldots,x_{i+1})$. Thus
we may assume $m$ is not divisible by $x_{i+1},\ldots,x_r$.
Since the monomial $mx_i$ belongs to $(I,x_r,\ldots ,x_{i+1})$, it must
belong to $I$. Since it has degree $e+1$, there must be a monomial
$m'\in I_e$ and an $l$ such that $mx_i=m'x_l$. Clearly $l\leq i$.
Since $m=(x_l/x_i)m'$, combinatorial stability of $I_e$ implies
that $m\in I_e$.
Thus 2a) of Theorem \ref{m-reg} holds as well, and $I$
is $e$-regular. \qed

{\bf Proof of Theorem 3:}
If $I$ is an ideal in generic coordinates which is $e$-regular,
then by Theorem \ref{m-reg}, $in_>(I)$ is generated in degrees $\leq e$,
where $<$ is the reverse lexicographic order.  By the above proposition,
$in(I_e)$ with respect to reverse lexicographic order is combinatorially
stable, and hence by Theorem \ref{char0-thm} we have:
\begin{eqnarray*}
\Delta(V_d(I)) & \leq & \delta(in_{>'}(V_d(gI))) \\
                    & \leq & max(2, \lceil reg(gI))/d \rceil) \\
                    & = & max(2, \lceil reg(I)/d \rceil)
\end{eqnarray*}
(where $g$ is a ``general" choice of coordinates, $>$ is reverse lexicographic
order, and $>'$ is the induced order),
proving Theorem \ref{reg-thm}.

\section{Comments on the main theorem}
\label{comments}

We have proved that for any
homogeneous ideal $I\subset S$, we have $\Delta (V_d(I))\leq max(2,\lceil
reg(I)/d \rceil )$. In particular, for suitable coordinates and order on $T_d$,
the Veronese ideal $V_d(I)$ has quadratic initial ideal for $d\geq reg(I)/2$,
and it follows that the Veronese subring $T_d/V_d(I)\subset S/I$  is Koszul
for $d\geq reg(I)/2$. In this section we will estimate the regularity of $I$
in order to bound $\Delta(V_d(I))$ in terms of other invariants of $I$
such as $\Delta (I)$. These
results can probably be improved, but we will give an example to show that
the most optimistic hopes are false.

\begin{thm}
\label{delta}
Let $r$ be the number of generators of the polynomial ring $S$.
 For any homogeneous ideal $I\subset S$,
$\Delta (V_d(I))\leq max(2,\lceil (r\Delta (I)-r+1)/d \rceil )$. In particular,
for suitable coordinates and order on $T_d$, the Veronese ideal $V_d(I)$
has quadratic initial ideal for $d\geq (r\Delta (I)-r+1)/2$.
\end{thm}

{\bf Proof:}
The Taylor resolution [Ta60] gives an upper bound on $reg(I)$,
specifically: $$reg(I) \leq r\Delta(I)-r+1.$$
With Theorem \ref{reg-thm}, this gives the result.
It is worth mentioning that, by Bayer and Stillman's Theorem \ref{m-reg},
there are actually upper and lower bounds relating $reg(I)$ and $\Delta (I)$:
$$\Delta (I) \leq reg(I) \leq r\Delta (I)-r+1.\mbox{    }\qed $$

The assumption $d\geq \Delta (I)/2$ is not enough to imply that
$V_d(I)$ has quadratic initial ideal, by the example at the end of this
section. We do not know the best estimate for $\Delta (V_d(I))$ in terms
of $\Delta (I)$. The problem is combinatorial in the sense that it suffices
to consider monomial ideals $I$.

We have been assuming that the field $k$ is infinite.
For arbitrary (in particular, finite) fields $k$, we have a slightly
weaker version of Theorem \ref{delta}: there is an order on $T_d$ such that
$V_d(I)$ has quadratic initial ideal for $d\geq r\lceil \Delta (I)/2 \rceil$.
We omit the proof, which is not too difficult given a definition of the
correct order.
The ordering which yields this result is defined as follows:

\begin{define}
For each monomial $m$ in $S$ of degree
$d$, we produce a vector
$$\nu(m) = (\nu_{11}(m),\nu_{12}(m),\ldots,\nu_{1r}(m),\nu_{21}(m),\ldots,
\nu_{2r}(m),\nu_{31}(m)\ldots),$$
where
$\nu_{ij}(m) = \cases{ 0 & if $x_j^i | m$ \cr 1 & else}$.
 The order on monomials in
$S$ of degree $d$ is then defined by $m > n$  if  $\nu(m) > \nu(n)$
in lexicographic order. We define the order of the variables in $T_d$ using
the above order on $S$.  Specifically, $z_m > z_n$ if $m > n$ in the order
on $S$ defined above.  Given this ordering on the variables in $T_d$, let the
order on the monomials in $T_d$ be reverse lexicographic order. \qed
\end{define}

For some monomial ideals $I$, we can improve the Taylor bound on the regularity
of $I$. First, since one direction of the proof of Proposition \ref{reg-stab}
does not use the Borel-fixed hypothesis, we have:

\begin{prop}
\label{stab-quick}
If $in(I)$ is generated in degrees $\leq u$ and $in(I)_u$ is combinatorially
stable, then $$reg(I)\leq u.$$
\end{prop}

Next we generalize the definition of combinatorial stability.

\begin{define}
Let $q$ be an integer.
A monomial ideal
$I$ is {\bf $q$-combinatorially stable}, if
for every $m \in I$ and
for each $j<max(m)$
there exists an integer $s$ with $1 \leq s \leq q$
such that $x_j^{s}/x_{max(m)}^{s}m \in I$.
\end{define}

\begin{prop}
Let $I$ be an ideal in generic coordinates, and let $e=$

\ni $\delta(in_>(I))$,
where $>$ is reverse lexicographic order.  If
$I$ is $q$-combinatorially stable,  then $reg(I) \leq e+(r-1)(q-1)$.
\end{prop}

{\bf Proof:}
Let $t = e+(r-1)(q-1)$.
By Proposition \ref{stab-quick}, we need only show that $J := in(I)_t$ is
 combinatorially stable.
Let $m \in J$.  $m = x_1^{b_1+c_1} \cdots x_r^{b_r+c_r}$, where $l :=
x_1^{b_1} \cdots x_r^{b_r} \in in(I)_e$, and set $n := x_1^{c_1}
\cdots x_r^{c_r}$. We have $\sum_{i=1}^r c_i = (r-1)(q-1)$.
If $max(m) = max(n)$, then $(x_i/x_{max(m)})m \in I_t$ for all $i$ with $1 \leq
i \leq x_{max(m)}$ because  $n$ is divisible by $x_{max(m)}$.
If $max(m) > max(n)$, then
$max(m) = max(l)$, and
either there exists some index $k$ such that $c_k \geq q$, or else
$c_j = q-1$ for all $j =
1,\ldots,n$.  In the first case, we can rewrite $m$ as $l'n'$, where
$l' = (x_k^s/x_{max(l)}^s)l$ and
$n' = (x_{max(l)}^s/x_k^s) n$, for
some $1 \leq s \leq q$.
After doing so, $max(m) = max(n)$ and we conclude as before.  In the second
case, the degree of $x_i$
in $x_i/x_{max(m)}m$
is $b_i + c_i+1$, and $c_i+1 = q$.  We may rewrite $(x_i/x_{max(m)})m$
as
 $(x_i^s/x_{max(l)}^s)ln'$, where $n'=
(x_{max(l)}^{s-1}/x_i^{s-1})n$, and $1 \leq s \leq q$ is chosen so that
$(x_i^s/x_{max(l)}^s)l \in I$. Thus
  $(x_i/x_{max(m)})m$ is in $I_t$.\qed

If $char \ k = 0$, then every ideal $I$ in generic coordinates is Borel-fixed
and hence 1-combinatorially stable. In this case, the proposition above
 yields $reg(I) \leq e$,
and in fact, equality holds, as Bayer and Stillman proved in [BaSt87].
In $char \ k = p$, every ideal $I$ in generic coordinates has a
 $q$-combinatorially stable
initial ideal for some $q$ that is a power of $p$
$\leq \delta(in(I))$ [Pa], but even if $q$ is chosen
to be as small as possible, $reg(I)$ can be strictly less than $e+(r-1)(q-1)$.
An example is the ideal $I = \{a^6,\,a^2b^4,\,a^2c^4,\,b^8,\,c^8\}$, which
is 8-combinatorially stable (implying a bound of 22 on the regularity), but
 has regularity 16.
Also, $I$ has a quadratic initial ideal in the Veronese embedding of degree 5,
which is strictly less than the degree of 7 given by Theorem \ref{char0-thm}.

As noted above, in
characteristic 0 and generic coordinates, the regularity of $I$ is
equal to $\delta(in(I))$, where the initial ideal is with respect to
reverse lexicographic order.
In characteristic $p$,
we cannot replace $reg(I)$ with $\delta(in(I))$ in the statement of Theorem
\ref{reg-thm}, as
the following example illustrates.

\begin{ex}
\label{counterexample}
2: A Borel-fixed ideal $I \subset k[a,b]$, with $char \ k= 2$,
 $e:= \delta(I) = 6$,
 such that the algebra $T_3/V_3(I)$ is not Koszul. It follows that
the initial ideal of the Veronese embedding of degree $\lceil e/2 \rceil = 3$
is not generated in degree 2
under any order and any generators for the graded algebra $T_3$. In fact
the ideal defined below has the same properties for $k$ of characteristic 0,
except that it is not Borel-fixed in characteristic 0.
\end{ex}

Let $I = (a^6, a^2b^4)$, and consider the embedding in degree $3 =
\lceil 6/2 \rceil$. Let $A=T_3/V_3(I)$. Thus, in the obvious coordinates
$y_i=a^{3-i}b^i$,
$$A=k[y_0,y_1,y_2,y_3]/(y_0^2=0,y_0y_2=y_1^2,y_0y_3=y_1y_2, y_1y_3=y_2^2=0).$$
The graded
vector space $Tor_3^A(k,k)$ is not entirely in degree 3: it has dimension
26 in degree 3 and dimension 2 in degree 4. So $A$ is not Koszul.

In fact, under the induced order used
throughout this paper, $in(V_3(I))$ requires 2 cubic generators.
However, $in(V_4(I))$ is generated in degree 2. The regularity of $I$
is 9.

\section{Resolution of multihomogeneous modules}
\label{resolution}

Fundamental to the discussion of rates above is the estimate of the rate
for a monomial ideal given (without proof) by Backelin in [Ba86].  The
case of quadratic monomials follows at once from the more precise
result of Fr\"oberg in [Fr75].  Fr\"oberg's result was recently reexamined
and reproved by Bruns, Herzog and Vetter [BrHeVe] using a different method.
J\"urgen
Herzog has pointed out to us that their method actually proves the
entire result claimed by Backelin, in a somewhat strengthened form,
and we now present this argument.

Let $S=k[x_1,\dots,x_r]$ be a polynomial ring over a field $k$.
We will regard $S$ as a $Z^r$-graded ring, graded by the
monomials.
Suppose that $I$ is
a monomial ideal of $S$, and set $A := S/I$;
the ring $A$ is again $Z^r$-graded.
If $M$ is a finitely generated
$Z^r$-graded module over $A$, then $M$ has a
$Z^r$-graded minimal free resolution over $A$.  The vector spaces
${\rm Tor}^A_i(k,M)$ are $Z^r$-graded.

For the purpose of bounding degrees it is convenient to
turn these multigradings into single gradings.  Rather than
simply using the total degree, we get a more refined result
by defining weights, as follows:
Let $w_1,\dots,w_r$ be non-negative real numbers.
For any monomial $m = x_1^{\alpha_1}\cdots x_r^{\alpha_r}$
define the {\it weight} of $m$ to be
 $w(m) =  \sum w_i\alpha_i$. Generalizing
the definition of $t^A_i(M)$ used above we define
$t^A_i(w, M)$ to be the maximal weight,
with respect to $w$, of a nonzero vector in
$Tor^A_i(k,M)$.

We can estimate the $t^A_i(w, M)$ as follows.
 Given an ordered set
$\{g_1,\ldots,g_s\}$ of generators of $M$ we get a
filtration
$$Ag_1 \subset Ag_1+Ag_2 \subset \cdots \subset Ag_1 + \cdots + Ag_s = M$$
of $M$ with quotients the cyclic modules
$A/J_i$ where
$J_i = ((g_1,\dots,g_{i-1}):g_i)$.
The set of generators $\{g_i\}$ also gives rise to a surjection $\phi$ of a
free
$A$-module $A^s$ to $M$ sending the $i^{\rm th}$ basis element to $g_i$.
It is easy to show that the kernel of $\phi$ has a filtration
whose successive quotients are
the ideals $J_i$
(see the proof of Theorem \ref{genBackthm}).
Thus the weights of the generators of the $J_i$
added to the weights of the $g_i$ give a bound for $t^A_1(w, M)$
(we get a bound and not an exact result because the set of
generators for the first syzygy of
$M$ produced from sets of generators for the $J_i$
may not be minimal).
Moreover, if we have a method for bounding the weights of syzygies
of the $J_i$, we may continue this process.
The following Lemma provides what we
require:

\begin{lem}
\label{ideallem}
Let $S=k[x_1,\dots,x_r]$ be a polynomial ring over a field $k$.
Suppose that $I$ is
an  ideal of $S$, generated by monomials $n_1,\dots,n_t$, and set $A := S/I$.

If $J = (m_1,\dots,m_s)\subset A$ is an ideal
generated by the images $m_i$ of monomials $m'_i$ of $S$, then the
quotient $((m_1,\dots,m_{s-1}):_Am_s)$ is generated by
the images in $A$ of
divisors of the monomials $m'_1,\ldots,m'_{s-1}$
and proper divisors of the
monomials $n_1,\ldots,n_t$.
\end{lem}

{\bf Proof:}
The quotient is the image in $A$ of
$((n_1,\dots,n_t,m_1,\dots,m_{s-1}):_S m_s)$ and is thus
generated by divisors of the monomials
$n_1,\dots,n_t,m_1,\dots,m_{s-1}$.  The divisors
of the $n_i$ that are not proper go to zero in $A$.\qed

Using Lemma \ref{ideallem} with the idea above we obtain:

\begin{thm}
\label{genBackthm}
Let $S=k[x_1,\dots,x_r]$ be a polynomial ring over a field $k$, and let
$w$ be a weight function on $S$ as above.  Suppose that $I$ is
an  ideal of $S$, generated by monomials, and set $A := S/I$.   Let $M$
be a $Z^r$-graded $A$-module
with  $Z^r$-homogeneous generators $\{g_1, \dots ,g_s\}$, of
weights $\leq d$, and set
$J_i = (Ag_1+\dots+Ag_{i-1}\,:\,_A\,g_i)$.
If the $J_i$ are generated by elements of weight $\leq e$, and both these
elements and the proper divisors of the generators of $I$ have
weights $\leq f$, then
for each integer $i \geq 1$ we have
$$
t^A_{i}(w, M) \leq d+e+(i-1)f.
$$
\end{thm}

{\bf Proof:}  We will inductively construct a
(not necessarily minimal)
free resolution
$$
\cdots \, \to F_2 \to F_1 \to F_0 \to M \to 0
$$
such that the generators of $F_0$ have weights $\leq d$,
the generators of $F_1$ have weights $\leq d+e$,
and for $i\geq 2$ the weights of the generators
of $F_i$ are $\leq d+e+(i-1)f$.
Since the minimal (multigraded) free resolution is a summand of
any free resolution, it follows that the weights of the
$i^{\rm th}$ free module in the minimal resolution are also
$\leq d+e+(i-1)f$, proving the desired inequality.

Let $F_0$ be $Z^r$-graded free $A$-module with
$s$ generators whose degrees match those
of the $g_i$, so that the surjection
$\phi_0:\, F_0\to M$
sending the $i^{\rm th}$ basis vector of $F_0$
to $g_i$
is multihomogeneous.  The weights
of the generators of $F_0$ are $\leq d$.

We will prove by induction on $s$ that the module
${\rm ker}(\phi_0)$ admits a filtration
with successive quotients isomorphic,
up to a shift in multidegree, to the ideals $J_i$, and that the
weight of the generators of this kernel are $\leq d+e$.
If we define $F'_0$ to be $F_0$ modulo the first
basis element, and define
$M'$ by the short exact sequence
$$
0 \to Ag_1 \to M \to M' \to 0
$$
then by the snake lemma we get a short exact sequence
$$
0 \to J_1 \to {\rm ker}(\phi_0) \to {\rm ker}( \phi'_0) \to 0,$$
 where $\phi'_0:\, F'_0 \to M'$ is the induced map.
By induction
${\rm ker}(\phi'_0) $ has a filtration with quotients
$J_2,\dots, J_s$, and generators of weights
$\leq d+e$.  This gives the desired filtration of
${\rm ker}(\phi_0) $.  The weights of
the generators of the copy of $J_1$ in the kernel
are  the weights of the monomials
generating the ideal $J_1$ plus the weight of $g_1$,
so they are also $\leq d+e$, and we are done.

Using this filtration of ${\rm ker}(\phi_0)$, we
define a free module $F_1$ whose generators have
weights $\leq d+e$
and a map $\phi_1:\, F_1\to F_0$ sending the generators
of $F_1$ to representatives $h_l$ in ${\rm ker}(\phi_0)$
of the generators of the successive quotients $J_i$.

We now repeat the argument, replacing $M$ by
${\rm ker}(\phi_0)$ and $\phi_0$ by $\phi_1$.
Lemma \ref{ideallem} applied to the
ideals $J_i$ implies that the
argument works as before if we replace $e$
by $f$: that is, ${\rm ker}(\phi_1)$
has a filtration with successive quotients
isomorphic (up to a shift in mult-degree)
to ideals with generators of weight $\leq f$.
This allows us to construct $F_2$ with generators of
degrees
$\leq d+e+f$ that map onto generators $h_i$ of
${\rm ker}(\phi_0)$ such that the ideals
$(Ah_1+\dots+Ag_{l-1}\,:\,_A\,g_l)$
have generators of weight $\leq f$.

We may continue to repeat the argument, using
the bound $f$ from the second step on, and
constructing the desired resolution. \qed

In the special case of the resolution of the residue class field $k$, we may
 take
$J_1$ to be the maximal ideal, and we get Backelin's result referred to above:

\begin{cor}
Let $S=k[x_1,\dots,x_r]$ be a polynomial ring over a field $k$.
Suppose that $I$ is
an  ideal of $S$, generated by monomials of degree
$\leq f$, and set $A := S/I$.  We have
$$
t^A_{i}(k) \leq 1+(i-1)(f-1).
$$
\end{cor}

Note that Theorem \ref{genBackthm} does {\em not}
prove much about the entries of the matrices
in even a non-minimal resolution.

With the ring $A$ as in the Theorem,
it would be interesting to know whether there is
 a minimal free resolution (say of the residue
class field ), with bases
for the free modules occurring, such that
the entries of the matrices representing the maps
of the resolution with respect to the given bases all have low
multidegree (or low weight).
The rate bound given in the above theorem, together with the requirement that
the resolution be minimal (so that each syzygy has weight at least 1) implies
that the individual entries of the $i^{th}$ syzygy matrix must have weights
bounded by $d+e+(i-1)f-i+1$.
We do not know whether there always exists
a free resolution with bases --- even non-minimal ---
such that the entries appearing in the matrices
are all proper divisors of the  generators of $I$, or even of
the least common multiple
of the generators of $I$.

\section{Segre Products of Veronese embeddings}
\label{multi}

The proofs of Proposition \ref{prop8} and Theorem \ref{char0-thm} can be
easily generalized to the Segre-Veronese case.  Below are the
definitions and statements we can make in this case.

Let $$S:=k[x_{11},\ldots,x_{1r_1},\ldots,x_{s1},\ldots,x_{sr_s}]$$ be the
coordinate ring of $\prs$, where
 $x_i = (x_{i1},\ldots,x_{ir_i})$ are the homogeneous coordinates
on $\P {r_i}$.
And let $$T:=k[\{z_{m}\}, m \mbox{ a monomial of $S$ of
multi-degree }(d_1,\ldots,d_s)]$$ be the coordinate ring of $\P N$.

\begin{define}
$\phi :T \longrightarrow S$
by $\phi(z_{m}) = m = m_1\cdots m_s$, where
$m_i = x_{i1}^{\alpha_{i1}} \cdots x_{is}^{\alpha_{is}}$.
\end{define}

${\rm ker}(\phi)$ is generated by the quadratic binomials of the form
$z_{m}z_{n}-z_{m'}z_{n'},$ where $m\cdot n =
m' \cdot n'$ in $S$. As in section \ref{vero},
if $a$, $b$ are monomials in $T$,
$a > b$ if $\phi(a) > \phi(b)$, or $\phi(a) = \phi(b)$ and $a>b$ in reverse
lexicographic order.

\begin{prop}
In reverse lexicographic order, the initial terms of the binomials
$z_{m}z_{n}-z_{m'}z_{n'}$
generate $in({\rm ker}(\phi))$.
\end{prop}

\begin{define}
The {\bf stabilization} $\{I\}$ of an ideal $I$ is defined to be
the ideal generated by
$$\{(x_j/x_{max(m)})m \ | m \in I, \ j=1,\ldots,max(m)\}.$$
\end{define}

\begin{define}
Call a multi-homogeneous monomial ideal $I$ {\it combinatorially stable} if
 it is
combinatorially stable in each set of variables independently.  That is, given
$$x_1^{\alpha_1}\cdots x_s^{\alpha_s} \in I,$$ we must have
$$\{x_1^{\alpha_1}\}\cdot \{x_2^{\alpha_2}\} \cdot \cdots \cdot
\{x_s^{\alpha_s}\} \subset I,$$ where $\{x_i^{\alpha_i}\}$ is the set of
all monomials necessary for an ideal in $k[x_{i0},\ldots,x_{ir_i}]$ containing
$x_i^{\alpha_i}$ to be combinatorially stable, and where the above product
 is the outer
product, i.e. all possible products of elements taken one from each set.
\end{define}

\begin{thm}
 If $I$ is a multi-homogeneous ideal whose initial ideal is
combinatorially stable, then $in(\sigma(I)) = \sigma(in(I))$ (where the
 initial terms on
the left are computed with respect to the induced order with ties broken by
reverse lex). Thus if $in(I)$ is generated in degrees $\leq
(u_1,\ldots,u_s)$ ($u_i$ = maximum degree of any generator with respect to
the $i^{th}$ set of variables), then $in(V(I))$ is generated in degrees $\leq
max(2,\lceil u_1/d_1 \rceil,\ldots,\lceil u_s/d_s \rceil)$.
\end{thm}

\section{ Another obstruction to having a quadratic initial ideal}
\label{obstsec}

In this section, we formulate a general obstruction to
the  existence of a quadratic initial ideal for a given polynomial ideal,
beyond the obvious requirement that the ideal must be generated by
quadratic polynomials, and even beyond the stronger requirement that
the quotient ring must be a Koszul algebra. We use the obstruction to show
that in certain dimensions, the ideal of
a generic complete intersection of quadrics
has no quadratic initial ideal, although every complete
intersection of quadrics is a Koszul algebra [BaFr85].

{\bf Note. }In discussing the existence of a quadratic initial ideal
for a homogeneous ideal $I$ in a polynomial ring $S=k[x_1,\ldots ,x_r]$,
we are asking whether there exists a set of coordinates $x_1',\ldots ,
x_r'$ (linear combinations of $x_1,\ldots ,x_r\in S_1$) and a monomial
order with respect to which $in(I)$ is generated by quadratic polynomials.

\begin{thm}
\label{obstruction}
Let $I$  be a homogeneous ideal in a polynomial ring $S$. Consider
the Krull dimensions $r=dim(S)$, $n=dim(S/I)$, $e=r-n$. (Thus, if $n\geq 1$,
$n$ is one more than the dimension of the projective variety defined by $I$.)
If there are coordinates and a monomial order such that the initial ideal
$in(I)$ has quadratic generators, then the ideal $I$ contains $e$ linearly
independent quadratic elements of the form:
\begin{eqnarray*}
c_1x_1^2 &=&x_2L_{1,2}+\cdots +x_{e+n}L_{1,e+n} \\
&\vdots \mbox{                                           }\\
c_{e}x_e^2 &=& x_{e+1}L_{e,e+1}+\cdots +x_{e+n}L_{e,e+n}
\end{eqnarray*}
for some basis $x_1,\ldots ,x_{e+n}$ for $S_1$ and some $c_i\in k$
and  linear forms $L_{ij}\in S_1$.

In particular, $I$ contains an $m$-dimensional space of quadrics of rank
$\leq 2(n+m)-1$ for every $m \leq codim(I)$.
\end{thm}

We recall that the rank of a quadratic form $Q$ over a field $k$
is the rank of a  symmetric matrix representing the form.

{\bf Proof. }We are given that there is a basis $x_1,\ldots ,x_{e+n}$
for the vector space $S_1$ and an ordering of the $x$-monomials, such that
the resulting initial ideal $in(I)$ is generated by $in(I)_2$. We can assume
that $x_1>\cdots >x_{e+n}$ in the monomial ordering.

We observe that for $i=1,\ldots ,e$,
there must be at least $i$ quadratic monomials $x_jx_k$
with $e-i+1 \leq j,k \leq e+n$ which are not allowable. Otherwise the Hilbert
series of $S/I$ would be at least equal to the Hilbert series
of an algebra
$k[x_{e-i+1},\ldots x_{e+n}]/(<i \mbox{ relations})$, so the dimension
of $S/I$ would be at least that of the latter ring, which is greater
than $n$; this contradicts $dim(S/I)=n$.

Thus, for $i=1,\ldots ,e$, there are $i$ monomials
$x_jx_k$, $e-i+1\leq j,k\leq e+n$, which are linear
combinations of earlier monomials $x_lx_m$. No matter what
monomial ordering we are using, at least one of $l$ and $m$
must be $> e-i+1$ in this situation. So, for all $1\leq i\leq e$,
 $R$ satisfies
$i$ independent relations of the form:
\begin{eqnarray*}
bx_{e-i+1}^2+(a_{e-i+2,1}x_{e-i+2}x_1+\cdots +a_{e-i+2,e+n}x_{e-i+2}x_{e+n})
  \cr
+(a_{e-i+3,1}x_{e-i+3}x_1+\cdots +a_{e-i+3,e+n}x_{e-i+3}x_{e+n})
+\cdots \quad \quad \cr
\quad \quad +(a_{e+n,1}x_{e+n}x_1+\cdots +a_{e+n,e+n}x_{e+n}^2)=0.
\end{eqnarray*}
This implies the statement of the lemma. \qed

\begin{cor}
Let $k$ be an infinite field, and let $I$ be an ideal in
$S=k[x_1,\ldots ,x_{e+n}]$ generated by $e$ generic quadratic forms defined
over $k$.
We assume that $n\geq 0$. (If $n\geq 1$, $S/I$
is the homogeneous coordinate ring of an $(n-1)$-dimensional
complete intersection of quadrics in {\it $\Proj ^{e+n-1}$}.)
If
$$n< \frac{(e-1)(e-2)}{6}$$
then generic complete intersection ideals $I$ as above do not admit any
quadratic initial ideal.
\end{cor}

{\bf Proof. }
Any $n$-dimensional complete intersection of homogeneous
quadrics in affine $(e+n)$-space can be described by a point of the
Grassmannian of $e$-dimensional subspaces of $S^2V$,
$\mbox{Gr}(X_e\subset S^2V)$, where we let $V=S_1$, a vector space of
dimension $e+n$ over $k$;
conversely, a nonempty open subset of this Grassmannian corresponds
to complete intersections.

Those $e$-dimensional linear spaces of quadrics which
generate an ideal which admits a quadratic initial ideal
can, by Lemma \ref{obstruction}, be written in the form:
\begin{eqnarray*}
c_1x_1^2 &=&x_2L_{1,2}+\cdots +x_{e+n}L_{1,e+n} \\
&\vdots & \quad \quad \quad \quad (*)\\
c_{e}x_e^2 &=& x_{e+1}L_{e,e+1}+\cdots +x_{e+n}L_{e,e+n}
\end{eqnarray*}
for some basis $x_1,\ldots, x_{e+n}$ for $V$ and some $c_i\in k$
and linear forms $L_{ij}\in V$. We want an upper bound for the
dimension of the space
of $e$-dimensional linear subspaces of $S^2V$ which can be written in this
form. If our bound is less than the dimension of the whole Grassmannian
$\mbox{Gr}(X_e\subset S^2V)$, then we will know that generic complete
intersections of quadrics in this dimension do not have any quadratic
initial ideal.

Our dimension estimate (which is not always sharp) is based on the following
observation. The basis $x_1,\ldots ,x_{e+n}$ for $V$ is not important,
only the flag $\langle x_{e+n},\ldots,x_{e+1}\rangle \subset \langle x_{e+n},
\ldots ,x_{e}\rangle \subset \cdots
\subset \langle x_{e+n},\ldots ,x_{1}\rangle =V$. That is, if a linear system
of quadrics
has the form $(*)$ for one basis $x_1,\ldots ,x_{e+n}$ of $V$, then it
has the same form for any basis which gives the same flag
$\langle x_{e+1},\ldots,x_{e+n}\rangle \subset\cdots $.

For any flag $V_n\subset V_{n+1}\subset \cdots \subset
V_{e+n}=V$, we consider an associated flag $W_{a_1}\subset W_{a_2}
\subset\cdots\subset W_{a_e}=S^2V$ defined by
$$W_{a_i}=(V_{n+i-1}\cdot V)+(V_{n+i}\cdot V_{n+i})\subset S^2V.$$
(That is, in terms of any basis $x_1,\ldots ,x_{e+n}$ adapted to the flag
$V_n\subset\cdots $, $W_{a_i}$ is the space of quadrics of the form
$c_{e-i+1}x_{e-i+1}^2=x_{e-i+2}L_{e-i+2}+\cdots +x_{e+n}L_{e+n}$, $L_j\in V$.)

{\bf Note. }We always use the notation $V_i$, $W_i$, etc.\ to denote
$i$-dimensional vector spaces.

Then, also associated to any flag $V_n\subset V_{n+1}\subset \cdots
\subset V_{e+n}=V$, we can consider the space of flags $X_1\subset
\cdots \subset X_e\subset V$ such that $X_i\subset W_{a_i}$ for
$i=1,\ldots ,e$.
Let $Q_{n,e}$ be the space of flags $V_n\subset V_{n+1}\subset \cdots
\subset V_{e+n}=V$ and $X_1\subset \cdots \subset X_e\subset S^2V$
such that $X_i\subset W_{a_i}$ for $i=1,\ldots,e$.
Then the image of the map
\begin{eqnarray*}
Q_{n,e}&\arrow &\mbox{Gr}(X_e\subset S^2V) \\
(V_i,X_i) &\mapsto &X_e
\end{eqnarray*}
contains the set of complete intersections of $e$ homogeneous quadrics
in $(e+n)$-space which have a quadratic initial ideal.

Moreover, it is easy to compute the dimension of $Q_{n,e}$, which is
an iterated projective bundle over the flag manifold $\mbox{Fl}
(V_n\subset V_{n+1}\subset\cdots\subset V_{e+n}=V)$: given a flag
$V_n\subset\cdots $ and hence the associated subspaces $W_{a_i}$,
we first choose the line $X_1\subset
 W_{a_1}$, then
a line $X_2/X_1\subset W_{a_2}/X_1$, and so on.

The dimension of the vector space $W_{a_i}$, for $i=1,\ldots,e$, is
\begin{eqnarray*}
a_i&= &(e+n)+(e+n-1)+\cdots +(e-i-2)+1 \\
 &= &  (e+n)(e+n+1)/2 -(e-i+1)(e-i+2)/2 +1.
\end{eqnarray*}
Using this, we compute the dimension of the variety $Q_{n,e}$:
\begin{eqnarray*}
\mbox{dim }Q_{n,e} &=& \mbox{dim Fl}(V_n\subset\cdots\subset V_{e+n}=V)
+\sum_{i=1}^e \mbox{dim }\Proj (W_{a_i}/X_{i-1}) \\
&=& en+\frac{e(e-1)}{2} +\sum_{i=1}^e[\frac{(e+n)(e+n+1)}{2}
-\frac{(e-i+1)(e-i+2)}{2}+1-i] \\
&=& en +\frac{e(e-1)}{2}+\frac{e(e+n)(e+n+1)}{2}-(\sum_{j=1}^e
\frac{j(j+1)}{2} )+e-\frac{e(e+1)}{2} \\
&=& e(n+\frac{(e+n)(e+n+1)}{2}-\frac{(e+1)(e+2)}{6}).
\end{eqnarray*}
So
\begin{eqnarray*}
& &\mbox{dim }Q_{n,e}<\mbox{dim Gr}(X_e\subset S^2V) \\
&\iff &e(n+(e+n)(e+n+1)/2-(e+1)(e+2)/6)<e((e+n)(e+n+1)/2-e) \\
&\iff &n<(e+1)(e+2)/6-e\\
&\iff &n<(e-1)(e-2)/6
\end{eqnarray*}

Thus, if $n<(e-1)(e-2)/6$, then the image of the map
$Q_{n,e}\arrow\mbox{Gr}(X_e\subset S^2V)$ has dimension less than
the dimension of $\mbox{Gr}(X_e\subset S^2V)$. Thus for
$n<(e-1)(e-2)/6$, the ideal generated by $e$ generic quadratic forms
in $e+n$ variables has no quadratic initial ideal. \qed

\vspace{.2in}
For example, the ideal generated by 3 generic
quadratic forms in 3 variables has no quadratic initial ideal.
Similarly for a generic complete intersection of 5 quadrics in
$\Proj ^5$, or a generic complete intersection of 6 quadrics in
$\Proj ^7$. (These last examples are smooth curves
of genus 129.)

Explicitly, say over a field $k$ of characteristic 0,
the ideal
$$I=(x(x+y),y(y+z),z(z+x))\subset k[x,y,z]=S$$
is a complete intersection of quadrics, and so $S/I$ is a Koszul algebra
[BaFr85], but one can check using Lemma \ref{obstruction}
that it has no quadratic initial ideal, for any coordinates
and any monomial order. (One has to check that no nonzero linear
combination of the relations $x(x+y)$, $y(y+z)$, $z(z+x)$ is
the square of a linear form, which is easy.)


\begin{thebibliography}{[DDDDDDD]}


\bibitem[An86]{} D. Anick:  On the homology of associative algebras,
{\em Trans. Am. Math. Soc.} {\bf 296} (1986), 641-659.

\bibitem[Ba86]{} J. Backelin:  On the rates of growth of the
homologies of Veronese
 subrings,  in {\em Algebra, Algebraic Topology, and Their Interactions}, ed.
 J.-E. Roos, Springer Lect. Notes in Math. {\bf 1183} (1986), 79-100.

\bibitem[BaFr85]{} J. Backelin and R. Fr\"{o}berg:  Koszul algebras, Veronese
subrings, and
 rings with linear resolutions, {\em Rev. Roum. Math. Pures  Appl.} {\bf 30}
 (1985),
85-97.

\bibitem[BaMa82]{}  S. Barcanescu and N. Manolache:  Betti numbers of
Segr\'{e}-Veronese
 singularities,  {\em Rev. Roum. Math. Pures  Appl.} {\bf 26} (1982) 549-565.

\bibitem[BaSt87]{} D. Bayer and M. Stillman: A theorem on refining division
orders by the reverse lexicographic order, {\em Duke Math. J.} {\bf 55}
(1987), 321-328.

\bibitem[BrHeVe]{} W. Bruns, J. Herzog, and U. Vetter: Syzygies and walks,
preprint.



\bibitem[Ei]{} D. Eisenbud:  Noncommutative Gr\"obner bases
for commutative algebras (to appear).

\bibitem[EiKoSt88]{} D. Eisenbud, J. H. Koh, and M. Stillman:  Determinantal
equations for curves
 of high degree, {\em Am. J. Math.} {\bf 110} (1988), 513-539.

\bibitem[ElKe90]{} S. Eliahou and M. Kervaire:
Minimal resolutions of some monomial
 ideals,  {\em J. Alg.} {\bf 129} (1990), 1-25.

\bibitem[Fr75]{}  R. Fr\"{o}berg:  Determination of a class of
Poincar\'{e} series. {\em Math. Scand.} {\bf 37} (1975), 29-39.

\bibitem[Ga79]{} A. Galligo:  Th\'eor\`eme de division et stabilit\'e, {\em
Ann. Inst. Fourier (Grenoble)} {\bf 24}(1979), 107-184.

\bibitem[Ke90]{} G. Kempf:  Some wonderful rings in algebraic geometry,
{\em J. Alg.} {\bf 134} (1990), 222-224.

\bibitem[La90]{}  R. Lazarsfeld:  Linear series on algebraic varieties,
 pp.\ 715-723 in
 {\em Proc. of the International Congress of Mathematicians, Kyoto,
 1990},  Springer-Verlag, New York.

\bibitem[Mu70]{}  D. Mumford:  Varieties defined by quadratic equations,
pp. 29-100 in
 {\em Questions on Algebraic Varieties}, Cremonese, Rome, 1970.

\bibitem[Pa]{} K. Pardue: Thesis, Brandeis University (in preparation).

\bibitem[Pa93]{} G. Pareschi:  Koszul algebras associated to adjunction
bundles, {\em J. Alg.} {\bf 157} (1993), 161-169.

\bibitem[Pr70]{}S. Priddy:  Koszul resolutions,
{\em Trans. Am. Math. Soc.} {\bf 152} (1970), 39-60.

\bibitem[Ta60]{} D. Taylor: Ideals generated by monomials in an $R$-sequence,
Thesis, University of Chicago, 1960.


\end{thebibliography}
\end{document}